\documentclass[12pt]{amsart}
\usepackage{amscd,amssymb}
\makeatletter
\@addtoreset{equation}{section}
\makeatother

\newcommand{\xref}[1]{{\rm \ref{#1}}}
\newcommand{\Supp}{\operatorname{Supp}}
\newcommand{\Exc}{\operatorname{Exc}}

\newcommand{\codim}{\operatorname{codim}}
\newcommand{\vv}{^{\operatorname{vert}}}
\newcommand{\hh}{^{\operatorname{hor}}}
\renewcommand{\emptyset}{\varnothing}
\newcommand{\down}[1]{\left\lfloor #1\right\rfloor}

\newcommand{\CC}{{\mathbb C}}
\newcommand{\QQ}{{\mathbb Q}}
\newcommand{\PP}{{\mathbb P}}
\newcommand{\ep}{\varepsilon}
\newcommand{\KKK}{\mathcal K}

\newtheorem{theorem}[equation]{Theorem}
\newtheorem{claim}[equation]{Claim}
\newtheorem{proposition}[equation]{Proposition}
\newtheorem{lemma}[equation]{Lemma}
\newtheorem{corollary}[equation]{Corollary}
\newtheorem{conjecture}[equation]{Conjecture}
\newtheorem{adjconjecture}[equation]{Adjunction Conjecture}

\theoremstyle{definition}
\newtheorem{example}[equation]{Example}
\newtheorem{definition}[equation]{Definition}
\theoremstyle{remark}
\newtheorem*{remark*}{Remark}

\author{Yuri G. Prokhorov}
\title{An application of the canonical bundle formula}
\address{
Department of Algebra, Faculty of Mathematics, Moscow State
Lomonosov University, Moscow 117234, Russia}
\address{
Max-Plank-Institut f\"ur Mathematik, Vivatsgasse 7, D-53111 Bonn,
Germany}
\thanks
{This work was is carried out  under the support of  grants RFFI
02-01-00441 and INTAS-OPEN 2000-269}

 \email{prokhoro@mech.math.msu.su}

\begin{document}
\begin{abstract}
We prove a part of Shokurov's conjecture on characterization of
toric varieties modulo the minimal model program and adjunction
conjecture.
\end{abstract}

\maketitle
\section{Introduction}

Let $X$ be a proper toric variety and let $D=\sum_{i=1}^n D_i$ be
the invariant divisor. It is well-known that
\begin{equation*}
n= \dim X+ \operatorname{rk} (\text{Weil divisors modulo numerical
equivalence}),
\end{equation*}
$K_X+D\sim 0$, and the pair $(X,D)$ has only log canonical
singularities.

V.V. Shokurov proposed that these properties can characterize
toric varieties:

\begin{conjecture}[{\cite{Sh1}}]
\label{conjecture}
Let $(X/Z\ni o,D=\sum d_iD_i)$ be a log variety such that $(X,D)$
has only log canonical singularities and $-(K_X+D)$ is nef over
$Z$. Then
\begin{equation}
\label{eq-conjecture-rho-lc-}
\sum d_i \le \sigma(X/Z)+\dim X,
\end{equation}
where $\sigma(X/Z)$ is the rank of the group Weil divisors on $X$
modulo the numerical equivalence over $Z$. Moreover, if the
equality holds, then $(X/Z\ni o,\down{D})$ is a toric log pair,
i.e, $(X/Z\ni o,\down{D})$ is formally (or analytically)
isomorphic to a toric log pair.
\end{conjecture}

\begin{example}
If $X/X\ni 0$ is a singularity germ, then $\sigma(X/X)$ is exactly
the rank of the group of Weil divisors modulo $\QQ$-Cartier
divisors (see \cite{Kaw}). For example, let $X/X\ni 0$ be the
singularity given by the equation $xy+zt=0$ in $\CC^4$ and let $D$
be the divisor cut out by $xy=0$. Then $(X,D)$ is log canonical,
$D$ has exactly four components and $\sigma(X/X)=1$. Thus we have
equality in \eqref{eq-conjecture-rho-lc-}. Clearly, the
singularity $X\ni 0$ is toric.
\end{example}

Conjecture \xref{eq-conjecture-rho-lc-} was proved by Shokurov in
dimension $2$ (see \cite{Sh1} and also \cite{Lect}). A special
case of this conjecture in dimension $3$ was verified in \cite{P}.
In this paper we prove a weak form of Shokurov's conjecture modulo
the log minimal model program (LMMP) and the adjunction conjecture
for fiber spaces. Before stating our main result we introduce
notation and definitions.

\subsection*{Notation}
We work over an algebraically closed field of characteristic zero.
Notation and conventions of the minimal model theory \cite{KMM},
\cite{Ut} will be used freely.

Let $(X, D)$ be a proper log pair. Write $D=\sum d_iD_i$ where the
$D_i$ are irreducible components. We denote
\begin{itemize}
\item[]
$\| D \|=\sum d_i$,
\item[]
$\langle D\rangle$ is the free abelian group generated by the
$D_i$,
\item[]
$\langle D\rangle^0$ is the subgroup of $\langle D\rangle$
consisting of $\QQ$-Cartier numerically trivial divisors,
\item[]
$
\sigma(X, D)= \operatorname{rk} \left(\langle D\rangle/ \langle
D\rangle^0\right).
$
\end{itemize}
Let $f\colon X\to Z$ be a fiber type contraction. We say that a
prime divisor $P$ is \emph{horizontal} if $f(P)=Z$ and
\emph{vertical} if $\dim f(P)<\dim Z$. For any divisor $D$ we have
the decomposition $D=D\hh+D\vv$ into the sum of \emph{horizontal}
and \emph{vertical parts}.

\begin{definition}
Let $(X,D)$ be a log pair. \emph{A numerical complement of
$K_X+D$} is a log divisor $K_X+D'$ with $D'\ge D$ such that
$K_X+D'$ has only log canonical singularities and numerically
trivial. We say that a log divisor $K_X+D$ is \emph{numerically
complementary} if it has at least one numerical complement.
\end{definition}

\begin{theorem}
\label{th-main-proj-ineq}
Let $(X, D)$ be a proper log variety such that $K_X+D$ is
numerically complementary. Assume that in dimension $\dim X$ both
the LMMP and Weak Adjunction Conjecture
\xref{conjecture-Adjunction} hold. Then
\begin{equation}
\label{eq-main-rho-lc}
\|D\|\le \sigma(X, D)+\dim X.
\end{equation}
Moreover, if the equality holds, then $X$ is rational. We can omit
the LMMP and Conjecture \xref{conjecture-Adjunction} as
assumptions when $\dim X \le 3$.
\end{theorem}

\begin{corollary}
Let $X$ be a non-singular proper variety of dimension $\le 3$ and
let $D$ be a reduced simple normal crossing divisor on $X$ such
that $K_X+D\equiv 0$. Then $\|D\|\le \rho(X)+\dim X$, where
$\rho(X)$ is the Picard number. Moreover, if the equality holds,
then $X$ is rational.
\end{corollary}

\begin{remark*}
According to recent preprint \cite{Am1} Weak Adjunction Conjecture
\xref{conjecture-Adjunction} hold in the case when $K_X+D$ has a
Kawamata log terminal numerical complement.
\end{remark*}

After finishing the main part of this work the author was informed
that similar and even more general results were obtained by J.
McKernan \cite{McK}. His proofs are completely different, much
easier, and mainly do not use the LMMP.

\subsection*{Acknowledgments}
The work has been completed during my stay at Max-Planck-Institut
f\"ur Mathematik in August 2002. I would like to thank MPIM for
hospitality and support. I also would like to thank Professor V.V.
Shokurov for some remarks and Dr. O. Fujino, who informed me about
\cite{McK}.

\section{Preliminary facts}
\subsection*{Adjunction (see \cite{K2}, \cite{Ambro}, \cite{F})}
Let $f\colon X\to Z$ be a contraction and let $D=\sum d_iD_i$ be a
$\QQ$-divisor on $X$ such that $d_i\le 1$ whenever $f(D_i)=Z$.

For a prime divisor $W\subset Z$, define a number $c_W$ as the log
canonical threshold:
\begin{equation*}
c_W=\sup\left\{c \mid (X, D+cf^*W)\ \text{is log canonical over
the generic point of $W$}\right\}.
\end{equation*}
Then the $\QQ$-divisor
\[
D_Z:=\sum_W (1-c_W)W
\]
is called the \emph{discriminant} of $K_X+D$.

Note that the definition of the discriminant $D_Z$ is a
codimension one construction, so computing $D_Z$ we can
systematically remove codimension two subvarieties in $Z$ and pass
to generic hyperplane sections $f_H\colon X\cap f^{-1}(H)\to Z\cap
H$. In particular, $f^*W$ is well-defined.

\begin{adjconjecture}[weak form]
\label{conjecture-Adjunction}
Let $f\colon X\to Z$ be a fiber space and let $D$ be a
$\QQ$-divisor on $X$ such that
\begin{enumerate}
\item
$(X, D)$ is log canonical,
\item
$K_X+D$ is $\QQ$-linearly trivial over $Z$.
\end{enumerate}
Then there is an effective $\QQ$-divisor $M_Z$ on $Z$ such that
$K_Z+D_Z+M_Z$ is log canonical and
\begin{equation}
\label{eq-canonical-bundle-formula}
K_X+D \mathbin{\sim_{\scriptscriptstyle{\QQ}}} f^*(K_Z+D_Z+M_Z).
\end{equation}
\end{adjconjecture}

According to Kawamata \cite{K2} this conjecture is true for
contractions of relative dimension one.

\subsection*{Properties of $\sigma$}

\begin{lemma}
\label{lemma-properties-sigma}
Let $g\colon X'\to X$ be a birational contraction, let $D'$ be the
proper transform of $D$, and let $E=\sum E_i$ be the exceptional
divisor. Then
\begin{equation}
\label{eq-vspomagat-varrho}
\sigma(X', D'+E)\le \sigma(X, D)+\| E\|.
\end{equation}
Moreover, if $g$ is a contraction of an extremal face in some
LMMP, then in \eqref{eq-vspomagat-varrho} the equality holds.
\end{lemma}
\begin{proof}
The surjective map $g_*\colon \langle D'+E\rangle\to \langle
D\rangle$ induces the exact sequence
\[
0\longrightarrow \langle E\rangle \longrightarrow \langle
D'+E\rangle/ g^*\langle D\rangle^0 \longrightarrow \langle
D\rangle/\langle D\rangle^0\longrightarrow 0
\]
Since $g^*\langle D\rangle^0\subset \langle D'+E\rangle^0$, this
gives us the desired inequality.
\end{proof}

\begin{corollary}
\label{corollary-sigma-extremal}
Let $g\colon X'\to X$ be a contraction of an extremal ray and let
$D$ be a divisor on $X'$. Then
\[
\sigma(X',D)-1\le \sigma(X,g_*D)=\sigma(X',D+E)-1\le \sigma(X',D).
\]
Furthermore, if $E$ is a component of $D$, then
$\sigma(X,g_*D)=\sigma(X',D)-1$.
\end{corollary}

\subsection*{Easy results in the local case}

\begin{proposition}[{\cite[18.22-23]{Ut}}]
\label{local}
Let $(X\ni o,B)$ be a log canonical singularity such that every
component of $B$ is $\QQ$-Cartier. Then
\[
\| B \| \le \dim X.
\]
Moreover, if the equality holds, then $X\ni o$ is a quotient of a
smooth point $X'$ by an abelian group $\mathfrak A$ acting on $X'$
free in codimension one. The proper transform $B_i'$ of each
component $B_i\subset \Supp(B)$ is $\mathfrak A$-stable.
Therefore, $(X,\down B)$ is a $\QQ$-factorial toric log pair.
\end{proposition}

\begin{corollary}
\label{corollary-local-codim}
Let $(X,B)$ be a log canonical log pair such that each component
$B_i$ of $B$ is $\QQ$-Cartier. Then $\|B\|\le \codim \cap B_i$.
\end{corollary}

\begin{corollary}[cf. {\cite[Corollary 18.24]{Ut}}]
\label{Fano-rho=1-Ut}
Let $V$ be a projective $\QQ$-factorial log terminal variety with
$\rho(V)=1$ and let $D$ be a boundary on $V$ such that $-(K_V+D)$
is nef and $(V,D)$ is log canonical. Then $\|D\|\le \dim V+1$.
Moreover, if the equality holds, then $(V,\down D)$ is a toric log
variety. If $V$ and all the $D_i$ are defined over a non-closed
field $\Bbbk$, then $V$ is $\Bbbk$-rational.
\end{corollary}

\begin{proof}
Assume that
\begin{equation}
\label{eq-lemma-D}
\|D\|\ge \dim V+1
\end{equation}
Take an embedding $V\subset \PP^n$ so that $V$ is projectively
normal and take a projective cone $X\subset \PP^{n+1}$ over $V$.
Let $B=\sum b_i B_i$ be the corresponding cone over $D$. Then $X$
is normal and $\QQ$-factorial. Let $\sigma\colon \tilde X\to X$ be
the blow up of the vertex, let $\tilde B$ be the proper transform
of $B$, and let $S$ be the exceptional divisor. We can write
\[
\sigma^*(K_X+B)=K_{\tilde X}+\tilde B-a(S,B)S.
\]
It is clear that $S$ is a Cartier divisor and $(S, \tilde
B|_S)\simeq (V,D)$. Since $K_V+D\equiv 0$, we have $a(S,B)=-1$.
Thus $\sigma^*(K_X+B)= K_{\tilde X}+\tilde B+S$ and by the
Inversion of Adjunction \cite[17.7]{Ut} the pair $(\tilde X,
\tilde B+S)$ is log canonical. Hence so is $(X,B)$. Moreover,
$(\tilde X, S)$ is plt (because $V$ is log terminal). Now take
$(X',B')$ and $\mathfrak A$ such as in Proposition \xref{local}.
Consider the diagram
\[
\begin{CD}
\tilde X @<{\tilde \pi}<< \tilde X'
\\
@V{\sigma}VV @V{\sigma'}VV
\\
X@<{\pi}<< X'
\end{CD}
\]
where $\tilde X'$ is the normalization of $\tilde X$ in the
function field of $X'$. Let $S'=\tilde\pi^{-1}(S)$. Then $\tilde
\pi$ is finite and $\tilde \pi|_{\tilde X'\setminus S'}$ is
\'etale in codimension $1$. By \cite[20.3]{Ut} the pair $(\tilde
X', S')$ is plt and $(\tilde X', \tilde B'+S')$ is log canonical,
where $\tilde B'=\tilde \pi^{-1}(\tilde B)$. In particular, $S'$
is irreducible and normal. Clearly each component of $\tilde B'$
is Cartier. Put $\Delta=\tilde B'|_{S'}$ and $\Delta_i=\tilde
B'_i|_{S'}$. Then $K_{S'}+\Delta$ is log canonical and numerically
trivial, $\Delta=\sum b_i \Delta_i$, $\sum b_i\ge \dim S'+1$, and
the $\Delta_i$ are ample numerically proportional Cartier
divisors. Hence $S'$ is a Fano variety of (Fano) index $\ge \dim
S'+1$ with only log terminal singularities. It is well-known (see
e.g. \cite[Th. 3.1.4]{IP}) that in this situation $S'\simeq
\PP^N$, $\sum b_i= \dim S'+1$, and all the $\Delta_i$ are
hyperplanes. This gives us the equality in \eqref{eq-lemma-D}.
Further, $V\simeq S=\PP^N/\mathfrak A$, where $\mathfrak A$ acts
on $S'=\PP^N$ so that all the $\Delta_i$ are stable. By Corollary
\xref{corollary-local-codim} we have $\cap \Delta_i=\emptyset$ and
$\down \Delta$ is a normal crossing divisor. This gives us that
$(V,\down D)$ is a toric pair. For the last statement, we note
that our construction is defined over $\Bbbk$. Since $\cap
\Delta_i=\emptyset$, we can take $\Bbbk$-coordinates on $\PP^N$ so
that the action of $\mathfrak A$ is monomial. The statement is
obvious in this case.
\end{proof}

\begin{lemma}
\label{lemma-local-extremal}
Let $f\colon X\to Z\ni o$ be the contraction of an extremal ray
$R$ and let $D=\sum d_iD_i$ be a boundary on $X$ such that all the
components $D_i$ are $\QQ$-Cartier and $-(K_X+D)$ is $f$-nef.
Assume that $X$ is $\QQ$-factorial and log terminal, and $(X,D)$
is log canonical. Put $D\hh=\sum_{D_i\cdot R>0} d_iD_i$ and
$D\vv=\sum_{D_i\cdot R\le 0} d_iD_i$. Then $\|D\vv\|\le \codim
f^{-1}(o)$. Moreover,
\begin{enumerate}
\item
if $f$ is divisorial, then $\|D\hh\|\le \codim f(\Exc(f))$;
\item
if $f$ is flipping and the flip $\chi\colon X\dashrightarrow
X^{{+}}$ exists, then $\|D\hh\|\le \codim f^{{+}{-1}}(o)$ and
$\|D\|\le \dim X+1$;
\item
if $f$ is of fiber type and if the LMMP holds in dimensions $<\dim
X$, then $\|D\hh\|\le \dim X/Z+1$. If furthermore the equality
holds, then a generic fiber $F$ is $\QQ$-factorial, $\rho(F)=1$,
the pair $(F,\down{D|_F})$ is toric, and $X_\KKK$ is
$\KKK$-rational, where $\KKK=\KKK(Z)$ is the function field of
$Z$.
\end{enumerate}
\end{lemma}
\begin{proof}
(i) follows by Corollary \xref{corollary-local-codim}.

(ii) Since $\chi (D\hh)\le \chi(D)\vv$, we have $\|D\hh\|\le
\codim f^{{+}{-1}}(o)$. Thus $\|D\|\le \codim f^{-1}(o)+\codim
f^{{+}{-1}}(o)$. On the other hand, by \cite[Lemma 5-1-17]{KMM} we
have $\codim \Exc(f)+\codim \Exc(f^+)\le \dim X+1$. This gives us
$\|D\|\le \dim X+1$.

(iii) Let $F$ be a generic fiber of $f$. Assume that $\|D\hh\|\ge
\dim X/Z+1=\dim F+ 1$. Denote $\Delta=D|_F$ and $\Delta_i=D_i|_F$.
Then all the $\Delta_i$ are ample and numerically proportional. We
claim that $F$ is $\QQ$-factorial and $\rho(F)=1$. Indeed, let
$g\colon \tilde F\to F$ be a small $\QQ$-factorialization and let
$\tilde \Delta=g^*\Delta$. Then $-(K_{\tilde F}+\tilde \Delta)$ is
nef, $(\tilde F,\tilde \Delta)$ is log canonical and $\tilde F$
has only log terminal singularities. Since $\tilde \Delta\neq 0$,
there is a $K_{\tilde F}$-negative extremal ray $R$. Then $\tilde
\Delta_i\cdot R>0$ for all $\tilde \Delta_i=g^*\Delta_i$ (here we
do not assume that $\tilde \Delta_i$ is irreducible). If
$\rho(\tilde F)>1$, this contradicts (i), (ii), or the inductive
hypothesis. By Corollary \xref{Fano-rho=1-Ut} we have
\[
\|D\hh\|= \|D|_F\|= \dim F+1
\]
and the pair $(F,\down{D|_F})$ is toric. Since $\|D\hh\|=
\|D|_F\|$, all the components of $D$ are defined over $\KKK$. So
the last assertion follows by Corollary \xref{Fano-rho=1-Ut}.
\end{proof}

\section{Proof of Theorem \xref{th-main-proj-ineq}}
Let $(X, D)$ be a log pair such that
\begin{equation}
\label{eq-main-contr}
\|D\|\ge \sigma(X, D)+\dim X.
\end{equation}
Assume that $K_X+D$ has a numerical complement $K_X+G$. Replace
$(X, G)$ with its minimal $\QQ$-factorial log terminal
modification (blow-up divisors with discrepancy $a(\cdot , G)=-1$)
and $D$ with the sum of its proper transform and the reduced
exceptional divisor. Thus $X$ is $\QQ$-factorial and log terminal.
Run $K_X$-MMP. By Lemma \xref{lemma-properties-sigma} and
Corollary \xref{corollary-sigma-extremal} this preserves
\eqref{eq-main-contr}. All the divisorial contractions are
positive with respect to $G$. Therefore, we cannot contract a
connected component of $G$ and $K\equiv -G$ cannot be nef. At the
end we get a fiber type extremal $G$-positive contraction $f\colon
X\to Z$. Note that our new $X$ is $\QQ$-factorial and log
terminal. It is sufficient to prove our theorem for this new $X$.

If $Z$ is a point, then $\rho(X)=1$ and $X$ is a log terminal Fano
variety . By Corollary \xref{Fano-rho=1-Ut}, $G=D$ and we have the
equality in \eqref{eq-main-contr}. Moreover, $(X,\down D)$ is a
toric pair.

Consider the case when $\dim Z>0$. On this step we use the
canonical bundle formula \eqref{eq-canonical-bundle-formula}.

\begin{proposition}
\label{prop-1-s-s}
Assumptions as in Theorem \xref{th-main-proj-ineq}. Assume
additionally that $X$ is $\QQ$-factorial, log terminal, and there
exists a fiber type extremal contraction $f\colon X\to Z$ with
$\dim Z>0$. Then the inequality \eqref{eq-main-rho-lc} holds.
Furthermore, if the equality holds, then $K_X+D$ is numerically
trivial over $Z$, a generic fiber $F$ is $\QQ$-factorial,
$\rho(F)=1$, the pair $(F,\down {D|_F})$ is toric, and both $X$
and $Z$ are rational.
\end{proposition}

\begin{proof}
Assume \eqref{eq-main-contr}. By Lemma \xref{lemma-local-extremal}
we have $\|D\hh\|\le \dim X/Z+1$. Since $f\colon X\to Z$ is an
extremal contraction, $D\vv=f^*\Delta$ for some effective
$\QQ$-Cartier divisor $\Delta$. Write
$D_Z=D_Z'+D_Z^{\prime\prime}$, where $D_Z'$ and
$D_Z^{\prime\prime}$ are effective divisors without common
components, $\Supp D_Z'\subset \Supp \Delta$ and none of the
components of $D_Z^{\prime\prime}$ is contained in $\Supp \Delta$.

\begin{claim}
\label{claim}
$\|D\vv\|\le \|D_Z'\|$.
\end{claim}
\begin{proof}
Let $W\subset Z$ be a prime divisor and let
$S=f^{-1}(W)_{\operatorname{red}}$. Since $\rho(X/Z)=1$, $S$ is
irreducible. Let $d$ be the coefficient of $S$ in $D\vv$. Write
$f^*W=kS$, where $k\in \mathbb N$. Then $d+c_Wk\le 1$ (because
$(X, D+c_Wf^*W)$ is log canonical over the generic point of $W$).
Hence, $d\le 1-c_W$. This proves the statement.
\end{proof}

Assume that $D\hh = 0$. Then $D\vv=D$ and
\[
\|D\vv\|\ge \sigma(X,D)+\dim X= \sigma(Z,\Delta)+\dim X.
\]
On the other hand, by \eqref{eq-canonical-bundle-formula} we have
that $K_Z+G_Z+M_Z$ is log canonical and numerically trivial. It is
clear that $D_Z'\le D_Z\le G_Z$, i.e., $K_Z+G_Z+M_Z$ is a
numerical complement of $K_Z+D_Z'$. Thus by the inductive
hypothesis and Claim \xref{claim},
\begin{equation*}
\| D\vv \|\le \| D_Z'\| \le \sigma(Z, D_Z')+\dim Z,
\end{equation*}
a contradiction.

Now assume that $D\hh \neq 0$. Then $\sigma(X,D) \ge
\sigma(X,D\vv)+1$ and $\|D\hh\|\le \dim X/Z+1$ (see Lemma
\xref{lemma-local-extremal}). Hence
\begin{multline}
\label{eq-to-surface-1}
\| D\vv \|\ge \sigma(X,D)+\dim Z-1\ge
\\
\sigma(X,D\vv)+\dim Z \ge \sigma(Z,\Delta)+\dim Z.
\end{multline}

As above we have
\begin{equation}
\label{eq-to-surface-2}
\| D\vv \|\le \| D_Z'\| \le \sigma(Z, D_Z')+\dim Z.
\end{equation}
This gives us the equalities in \eqref{eq-to-surface-1},
\eqref{eq-to-surface-2} and \eqref{eq-main-contr}. In particular,
$ \| D_Z'\| = \sigma(Z, D_Z')+\dim Z$. Since $D_Z'$ is numerically
complementary, the inductive hypothesis give us that $Z$ is
rational.

Finally, by Lemma \xref{lemma-local-extremal} the pair $(F,
\down{D|_F})$ is toric and $X$ is rational.
\end{proof}

Since by \cite{K2} Conjecture \xref{conjecture-Adjunction} holds
when $\dim X/Z=1$, to prove the last part of the theorem we have
to consider only the case when $\dim X=3$ and $Z$ is a curve. In
this case, Theorem \xref{th-main-proj-ineq} is an immediate
consequence of the following.

\begin{proposition}
Let $f\colon X\to Z$ be an extremal contraction to a curve, where
$X$ is $\QQ$-factorial and log terminal, and let $D$ be a boundary
on $X$ such that $K_X+D$ is numerically complementary. Then $\|
D\|\le \dim X+\sigma(X,D)$. Furthermore, if the equality holds,
then $\sigma(X,D)=2$, $K_X+D\equiv 0$, and $X$ is rational.
\end{proposition}
\begin{proof}
Let $K_X+G$ be a numerical complement. Assume that
\begin{equation}
\label{eq-to-curve-contr}
\|D\|\ge \sigma(X, D)+\dim X.
\end{equation}
Since $Z$ is a curve, $\rho(X)=2$ and $\sigma(X,D)\le 2$. By Lemma
\xref{lemma-local-extremal},
\begin{equation}
\label{eq-rho=2-Dhh}
\|D\hh\|\le \| G\hh \|\le \dim X, \qquad \|D\vv\|\ge \sigma(X, D).
\end{equation}
As $\rho(X)=2$, the Mori cone $\overline{NE}(X)$ is generated by
two extremal rays, say $R$ and $Q$. Let $R$ is the ray
corresponding to the curves in fibers of $f$.

\begin{claim}
There is a boundary $\Delta\le G$ such that $(X,\Delta)$ is a
Kawamata log terminal log Fano variety.
\end{claim}

\begin{proof}
Obviously, $G\vv\cdot Q>0$. Thus $(K+G-\ep G\vv)\cdot Q<0$ for
$\ep>0$. Since $(K+G-\ep G\vv)\cdot R=0$, for $0<\ep_1\ll \ep$, we
have $(K+G-\ep G\vv-\ep_1G\hh)\cdot R<0$ and $(K+G-\ep
G\vv-\ep_1G\hh)\cdot Q<0$. Put $\Delta=G-\ep G\vv-\ep_1G\hh$. Then
$K+\Delta$ is anti-ample and Kawamata log terminal.
\end{proof}

By the Cone Theorem the ray $Q$ is contractible (i.e., there is a
$D\vv$-positive extremal contraction $g\colon X\to W$). This
implies that $Z\simeq\PP^1$. It is clear that $g$ cannot have
fibers of dimension $\ge 2$. Since all the components of $G\vv$
are strictly positive with respect to $Q$, by Lemma
\xref{lemma-local-extremal} we have
\begin{equation}
\label{eq-rho=2-Dver}
\|D\vv\|\le \|G\vv\|\le 2
\end{equation}

If $\sigma(X,D)=1$, then all the components of $D$ are numerically
proportional. By \eqref{eq-rho=2-Dhh} and
\eqref{eq-to-curve-contr}, $D\neq D\hh$. Hence $D=D\vv$. This
contradicts \eqref{eq-rho=2-Dver}. Thus $\sigma(X,D)=2$ and
$\|D\|\ge \dim X+2$. Combining \eqref{eq-rho=2-Dhh},
\eqref{eq-rho=2-Dver} and \eqref{eq-to-curve-contr} we get
$\|D\vv\|= 2$, $\|D\hh\|= \dim X$, $\|D\|= \dim X+2$, $D=G$ and
$K_X+D\equiv 0$. The rest follows by Lemma
\xref{lemma-local-extremal}.
\end{proof}

In conclusion, note that our method allow us to prove Conjecture
\xref{conjecture} in the general form (at least in dimension
three). We do not include these results because now there is much
better approach \cite{McK}.

\end{document}